\documentclass[12pt, leqno]{article}

\usepackage{amssymb}
\input epsf
\def\q \m#1#2{{\raise1pt\hbox{$#1$}\kern-1pt\big/
               \kern-1pt\raise-1pt\hbox{$#2$}}}

\def\bN{{\rm \bf N}}
\def\bR{{\rm \bf R}}
\def\bZ{{\rm \bf Z}}
\def\bQ{{\rm \bf Q}}
\def\bC{{\rm {\bf C}}}
\def\ch{{\rm  ch}}
\def\bH{{\rm \bf H}}
\def\sR{{ \rm \scriptsize  \bf R}}
\def\sZ{{ \rm \scriptsize  \bf Z}}

\newfam\msbfam

\font\twelmsb=msbm10 at 12pt
\font\tenmsb=msbm10 at 10 pt
\font\sevenmsb=msbm10 at 7pt
\font\fivemsb=msbm10 at 5pt

\textfont\msbfam=\twelmsb
\textfont\msbfam=\tenmsb
\scriptfont\msbfam=\sevenmsb
\scriptscriptfont\msbfam=\fivemsb\newtheorem{thm}{Theorem}[section]

\newtheorem{lemma}{Lemma}[section]

\newtheorem{cor}{Corollary}[section]

\newcommand{\g}{{\frak g}}

\begin{document}

\renewcommand{\theequation}{\thesection.\arabic{equation}}
\setcounter{equation}{0}

\centerline{\Large {\bf On Family Rigidity Theorems II}}
\vskip 10mm  
\centerline{\bf Kefeng   LIU and Xiaonan MA}
\vskip 8mm

{\bf Abstract.} In [{\bf LM}], we proved a family version of the famous Witten
 rigidity theorems and several family vanishing theorems for elliptic genera. 
In this paper, we gerenalize our theorems [{\bf LM}] in two directions. 
First we establish a family 
rigidity theorem for the Dirac operator  on loop space twisted by general 
positive energy loop group representations. Second we prove a family 
rigidity theorem for $spin^c$-manifolds. 
Several vanishing theorems on both cases are also obtained.\\

{\bf 0 Introduction } In [{\bf W}], Witten considered the indices of elliptic 
operators on the free loop space ${\cal L} M$ of a manifold $M$. 
In particular the index of  the formal signature operator on loop space is 
exactly the elliptic genus of Landweber-Stong.
 Witten made the conjecture about the rigidity of these elliptic 
operators which says that 
their $S^1$-equivariant indices on $M$ are independent of $g\in S^1$.  
 We refer the reader to [{\bf T}], [{\bf BT}], [{\bf H}], [{\bf K}], [{\bf L}] and [{\bf O}] for the history of the subject.

In [{\bf Liu2}], the first author observed that these rigidity theorems are 
consequence of their modular invariance. This allowed him to give a simple 
and unified proof of the above conjectures of Witten. 
In [{\bf Liu4}], it was proved the rigidity of 
 the Dirac operator on loop space twisted by positive energy loop group 
representations of any level, while the Witten rigidity theorems are 
the special cases of level 1.  An $\widehat{\frak U}$-vanishing 
theorem for loop spaces with spin structure, which is an analogue of the 
famous $\widehat{\frak U}$-vanishing theorem of Atiyah and Hirzebruch
[{\bf AH}], was also proved in [{\bf Liu4}]. 
Recently, by using Liu's idea, Dessai [{\bf D1}] proved 
a version of rigidity theorem for $spin^c$-manifolds.

The purpose of our paper is to generalize these results to family case.

 Let $M,B$ be two compact smooth manifolds, 
and $\pi: M\to B$ be a  submersion with compact fibre $X$.
Let a compact Lie group $G$ act fiberwisely on $M$, that is the action 
preserves each fiber of $\pi$. Let $P$ be a family of elliptic operators
along the fiber $X$, commuting with the action $G$. Then the family index 
of $P$ is 
\begin{eqnarray} 
{\rm Ind} (P) = {\rm Ker } P - {\rm Coker } P \in K_{G} (B).
\end{eqnarray}

Note that ${\rm Ind} (P)$ is a virtual $G$-representation. 
Let $\ch_g ({\rm Ind} (P))$ with $g\in G$ be the equivariant 
Chern character of ${\rm Ind} (P)$ evaluated at $g$.

To consider rigidity, we only need to restrict to the case when $G=S^1$. 
>From now on we let $G=S^1$. A family elliptic operator $P$ is called 
{\em rigid on equivariant Chern character level} 
with respect to this $S^1$-action, if $\ch_g ({\rm Ind}(P)) \in H^* (B)$
 is independent of $g \in S^1$.

In [{\bf LM}], several family rigidity and vanishing results for elliptic 
genera were obtained. As pointed out in [{\bf LM}], by taking expansions 
in $H^*(B)$, from the family rigidity and vanishing theorems we get  many
 higher level rigidity and vanishing results for characteristic numbers of
 the family. These characteristic numbers may not be the indices of any 
elliptic operators.

 This paper is the continuation of [{\bf LM}], 
and is naturally divided into two parts.
In Section 1, we prove a family rigidity theorem of the Dirac operator 
on loop space twisted by positive energy loop group representations, 
and we also derive some vanishing theorems.
In Section 2, we prove the family rigidity and vanishing theorems 
for $spin^c$-manifolds which generalize the results of Dessai [{\bf D1}]

. 

\newpage

\section{ \normalsize Loop groups and family rigidity theorems}

\setcounter{equation}{0}

This Section is organized as follows: 
In Section 1.1, we recall the modular invariance of the characters 
of the representations of affine Lie algebra. 
 In Section 1.2, we state  the family rigidity theorems of the Dirac operator 
on loop space twisted by general positive energy loop group representations 
for  spin manifold.
In Section 1.3,  we prove the main theorem, Theorem 1.2. 
In Section 1.4, we derive some vanishing theorems.

\subsection{ \normalsize  Characters of affine Lie algebras}

Let $G$ be a simple, simply connected compact Lie group and $LG$ be 
its loop group. There is a central extension $\widetilde{L} G$ of $G$
\begin{eqnarray}
1 \to S^1 \to \widetilde{L} G \to LG \to 1.
\end{eqnarray}
The circle group $S^1$ acts on $LG$ by the rotation $R_{\theta}$, 
$R_{\theta} \nu (\theta') = \nu (\theta'-\theta)$. The action of $S^1$ on 
$LG$ lifts (essentially uniquely) to an action on $\widetilde{L} G$.
We say the representation $U$ of $LG$ is symmetric if 
$R_{\theta}U_{\nu}R_{\theta}^{-1}= U_{R_{\theta} \nu}$. 
We say a representation $E$ of $\widetilde{L} G$ is positive energy 
[{\bf PS}, Chap 9] if

(a) $E$ is a direct sum of irreducible representations;

(b) $E$ is symmetric and $E^0 = \oplus_{j \in \bN} E_j$ is dense in $E$,
 where $E_j= \{v\in E :  R_{\theta}v = e^{-ij v}\}$ and $E_j$ is a finite
 dimensional complex representation of $G$;

(c) The action of $\widetilde{L} G \rtimes S^1$ on $E$ naturally extends
 to a smooth action of $\widetilde{L} G \rtimes {\rm Diff}^+(S^1)$, where 
$ {\rm Diff}^+(S^1)$ is the group of orientation preserving diffeomorphisms 
of $S^1$.\\

Let $\g$ be the Lie algebra of $G$. Let $\eta$ be the Cartan subalgebra, 
$W$ be the  Weyl group of $\g$.  Denote by $Q= \Sigma_{i=1}^l \bZ \alpha_i$, 
where $\{ \alpha_i\}$ is the root basis, the root lattice of $\g$.
 Then the affine Lie algebra associated to $\g$ is
\begin{eqnarray}
\widehat{L}  \g = \g \otimes_{\bR} \bC [t, t^{-1}] \oplus \bC K 
\oplus \bC  d,
\end{eqnarray}
where $K$ (resp. $d$)  is the infinitesimal generator of the central 
element (resp. the rotation of $S^1$) of $\widetilde{L} G$. 
$\widehat{L}  \g $ has the triangle decomposition 
\begin{eqnarray}
\widehat{L}  \g = \widehat{\eta}_- \oplus \widehat{\eta}
\oplus \widehat{\eta}_+,
\end{eqnarray}
where $\widehat{\eta}_{\pm}$ are the nilpotent subalgebras and $\widehat{\eta}
= \eta\otimes_{\bR} \bC \oplus \bC K \oplus \bC  d$ is the Cartan subalgebra. 
Let $(, )$ be the normalized symmetric invariant bilinear form on
 $\widehat{L}  \g $ which extends the standard  symmetric bilinear 
form  on $\g$, such that
\begin{eqnarray}
 \qquad (\bC K \oplus \bC d, \g \otimes_{\bR} \bC [t, t^{-1}]) =0; 
\quad (K, K)=0; \quad (d,d)=0; \quad (K, d)=1.
\end{eqnarray}
Let $\widehat{\eta}^*$ be the dual of $\widehat{\eta}$ with respect to $(, )$. 
Let $\left \langle , \right \rangle$ denote the pairing between 
$\widehat{\eta}^*$ and $\widehat{\eta}$. 
Then the level of $\lambda \in \widehat{\eta}^*$
is defined to be $ \left \langle \lambda, K \right \rangle$.
Let $\Lambda_0, \delta \in \widehat{\eta}^*$ be the elements such that
$\delta|_{\eta\oplus \bC K} =0, \left \langle \delta, d \right \rangle=1$;
$\Lambda_0|_{\eta\oplus \bC d} =0, \left \langle \Lambda_0, K\right \rangle=1$.

It is known that $\widehat{L}  \g $ falls into class $X^{(1)}_N$ in the classification of Kac-Moody algebras [{\bf Kac}, Chap 7]. 
An $\widehat{L}  \g $-module $V$ is called a highest weight module with 
the highest weight $\Lambda\in \widehat{\eta}^*$ if there exists a non-zero
 vector $v\in V$ such that
\begin{eqnarray}
\widehat{\eta}_+ (v) =0; \quad  h(v) = \Lambda (h) v, \quad 
{\rm  for  }\quad  h\in \widehat{\eta};  
\quad {\rm  and   }\quad  U(\widehat{L}  \g ) (v) =V.
\end{eqnarray}
where $U(\widehat{L}  \g ) $ is the universal envelopping algebra of 
$\widehat{L}  \g $. An irreducible representation $L(\Lambda)$ of $
\widehat{L}  \g $ with the highest  weight $\Lambda$ is called of level $k= 
\left \langle \Lambda, K \right \rangle$. $L(\Lambda)$ is said to be 
integrable if $\Lambda \in P_+ =
\{ \lambda \in \widehat{\eta}^*: \left \langle \lambda, \alpha_i \right \rangle
\in \bN \quad  {\rm for \,\,  all }\quad i\} $,  the set  of dominant integral
 weights. 
An integrable highest weight representation $L(\Lambda)$ of $\widehat{L} \g$
 can always be lifted to a representation of $\widetilde{L} G$ which turns out
 to be  irreducible and of positive energy.

Since each $ \widetilde{L} G$-module $V$ has a weight space decomposition 
$V= \oplus_{\lambda \in \widehat{\eta}^*} V_{\lambda}$, we can define 
formal Kac-Weyl character of $U$ as 
$\ch_V= \Sigma_{\lambda \in \widehat{\eta}^*} (\dim(V_{\lambda}) ) e^\lambda$.

The normalized character of $L(\Lambda)$ is $\chi_\Lambda = q ^{m_\Lambda}
\ch_{L(\Lambda)}$, where [{\bf Kac}, (12.8.12)]
\begin{eqnarray}
m_{\Lambda} = { (\Lambda + 2 \rho, \Lambda) \over 2 (m + h^\vee)}
- { m \dim \g \over 24 (m + h^\vee)},
\end{eqnarray}
where $h^\vee = \left \langle \rho, K \right \rangle$, the dual Coxeter number,
  $\rho= \overline{\rho} + h^\vee \Lambda_0$ [{\bf Kac}, (6.2.8)], 
and $\overline{\rho}$ is  half the sum of the positive roots of $\g$. 
We call $q^{m_\Lambda}$ the anomaly factor.

Let $M= \bZ(W\cdot \theta)$ be a lattice in $\eta^*$, where $\theta$ is the 
long root in $\eta$, and $W$ is the Weyl group of $\g$. For any integer $m$,
 let $P^m_+= \{\lambda \in P_+| \left \langle \lambda, K \right \rangle = m\}$
 be the level $m$ subset of the dominant integral weights.

If we choose an orthonormal basis $\{ v_j \}_{j=1}^l$ of $\eta^* \otimes_{\sR}
 \bC$, such that for $v \in \widehat{\eta}^*$, then we have
$$v= 2 \pi i ( \Sigma_{j=1}^l z_j v_j - \tau \Lambda_0 + u \delta).$$
we denote $z= \Sigma_{s=1}^l z_s v_s \in  \eta^* \otimes_{\sR}
 \bC$. Recall the classical theta functions associated to the lattice $M$ 
is defined by
\begin{eqnarray}
\Theta_{\lambda} (z, \tau) = e^{ 2 \pi i m u} 
\sum_{\gamma \in M + m^{-1} \overline{\lambda}} 
e^{\pi i m \tau (\gamma, \gamma) 
+ 2 \pi i  m (\gamma, z)}.
\end{eqnarray}
Here $\overline{\lambda}$ means the orthogonal projection of $\lambda $ from 
$\widehat{\eta}^*$ to $ \eta^* \otimes_{\sR} \bC$ with respect to 
the bilinear form $(\cdot, \cdot)$, and $\gamma = \Sigma_{i=1}^l \gamma_i v_i$
 with $(\gamma, z) = \Sigma_{i=1}^l \gamma_i z_i$. Then we can  express
$\chi_{\Lambda}$ as a finite sum 
\begin{eqnarray}
\chi_{\Lambda} (z, \tau) = \sum_{\lambda \in P^m {\rm \scriptsize   mod} 
(m M + {\scriptsize  \rm \bf C} \delta)} c^{\Lambda}_{\lambda}(\tau) 
\Theta_{\lambda} (z, \tau),
\end{eqnarray}
Where $P^m$ is the level $m$ element in the integral weight lattice, 
and $\{c^{\Lambda}_{\lambda}(\tau)\}$ are some modular forms of weight 
${- {1 \over 2} l}$, which are called string functions in 
[{\bf Kac}, \S 12.7, \S 13.10].

One of the important facts about the formal character is that 
$\chi_{\Lambda}$ is holomorphic in $\{ v\in \widehat{\eta}: 
{\rm Re} (\delta, v)>0\}= \{(z,\tau, u)\in \bC^{l+2}, {\rm Im}(\tau) >0 \}$.

Now, we state the following important Kac-Peterson theorem on the modular
 transformation property of $\chi_{\Lambda}$ under $SL_2(\bZ)$ 
[{\bf Kac}, Theorem 13.8].
\begin{thm} Let $\Lambda \in P^m_+$. Then
\begin{eqnarray}
\chi_{\Lambda} ({z\over  \tau}, -{1 \over \tau}) = e^{\pi i m (z,z) /\tau}
 \sum_{\Lambda' \in P^m_+ {\rm  \scriptsize mod } 
{\rm  \scriptsize \bf C}\delta} S_{\Lambda, \Lambda'} 
\chi_{\Lambda'} (z, \tau),
\end{eqnarray}
for some complex numbers $S_{\Lambda, \Lambda'}$, and 
\begin{eqnarray}
\chi_{\Lambda} (z, \tau + 1) = e^{2 \pi i m_{\Lambda}} 
\chi_{\Lambda} (z, \tau).
\end{eqnarray}
\end{thm}

 By (1.8), for $\alpha \in M$, we also have 
\begin{eqnarray}  \begin{array}{l}
\chi_{\Lambda} (z+ \alpha, \tau) = \chi_{\Lambda} (z, \tau);\\
\displaystyle{
\chi_{\Lambda} (z+ \alpha \tau, \tau) = e^{ 2 \pi i m (z, \alpha)
+ \pi i m (\alpha, \alpha)} \chi_{\Lambda} (z, \tau).}
\end{array}\end{eqnarray}
This, together with its transformation formulas (1.9), (1.10), means that 
$\chi_{\Lambda}$ is an $l$-variable Jacobi form of index $m/2$ and weight $0$.

\subsection{ \normalsize Family rigidity theorem of general elliptic genera}

Let $\pi: M\to B$ be a  fibration of compact manifolds with fiber $X$, 
and $\dim X= 2k$. We assume that the $S^1$ acts fiberwisely on $M$,
and $TX$ has an $S^1$-equivariant spin structure.
Let $\Delta(TX) = \Delta^+ (TX) \oplus \Delta ^- (TX) $ be 
the spinor bundle of $TX$.
Let $D^X$ be the Dirac operator on $\Delta(TX)$ 
which is defined fiberwisely on the fiber $X$. 

For a vector bundle $F$ on $M$, we let
\begin{eqnarray}\begin{array}{l}
S_t (F) = 1 + t F + t^2 S^2 F + \cdots,\\
\Lambda_t (F) = 1 + tF + t^2 \Lambda^2 F + \cdots,
\end{array}\end{eqnarray}
be the symmetric and respectively the exterior power  operations
 in $K(M)[[t]]$.

Assume $E$ is an irreducible positive energy representation of 
$\widetilde{L} Spin(2l)$ and $V$ is an $S^1$ equivariant vector bundle with structure group $Spin(2l)$ over $M$. 
Let $\Lambda$ be the highest weight of $E$ and $m$ be the level of $E$. 
 By the discussion in Section 1.1, we have the decomposition 
$E= \oplus_{n\geq 0} E_n$ under the action of $R_{\theta}$. 
Here $E_n$ is a finite dimensional representation of $Spin(2l)$. Let $P$ be the frame bundle of $V$, which is a $Spin(2l)$ principal bundle. We define
\begin{eqnarray}
\psi (E, V) = \Sigma_{n\geq 0} (P\times_ {Spin(2l)} E_n) q^n \in K(M)[[q]].
\end{eqnarray}

 Let $p_1(\cdot)_{S^1}$ 
denote the first $S^1$-equivariant Pontrjagin class.

\begin{thm} For $E$ an irreducible positive energy representation of 
$\widetilde{L} Spin(2l)$ of highest weight of level $m$, 
if $p_1(TX)_{S^1} = m p_1(V)_{S^1}$, then  the elliptic operator 
$$D^X \otimes _{m=1}^\infty S_{q^m} (TX) \otimes \psi (E, V) $$
is rigid on equivariant Chern character level.
\end{thm}

Theorem 1.2 actually holds for any semi-simple and simply connected Lie group,
instead of $Spin(2l)$.

\subsection{\normalsize Proof of Theorem 1.2}

For $\tau \in \bH = \{ \tau \in \bC; {\rm Im} \tau >0\}$,
 $q= e^{ 2\pi i \tau}$, $v\in \bC$,  let  
\begin{eqnarray}\begin{array}{l}
\theta_3(v, \tau)=c(q)\Pi_{n=1}^\infty (1 + q^{n-1/2} e^{2 \pi i v}) 
\Pi_{n=1}^\infty (1 + q^{n-1/2} e ^{-2 \pi i v}),\\
\theta_2(v, \tau)=c(q)\Pi_{n=1}^\infty (1 - q^{n-1/2} e^{2 \pi i v}) 
\Pi_{n=1}^\infty (1 - q^{n-1/2} e ^{-2 \pi i v}),\\
\theta_1(v, \tau)=c(q) q^{1/8} 2 \cos(\pi  v)
\Pi_{n=1}^\infty (1 + q^{n} e^{2 \pi i v}) 
\Pi_{n=1}^\infty (1 + q^{n} e ^{-2 \pi i v}),\\
\theta (v, \tau)=c(q)q^{1/8} 2 \sin (\pi v )
\Pi_{n=1}^\infty (1 - q^{n} e^{2 \pi i v}) 
\Pi_{n=1}^\infty (1 - q^{n} e ^{-2 \pi i v}).
\end{array}\end{eqnarray}
be the classical Jacobi theta functions [{\bf Ch}], 
where $c(q)= \Pi_{n=1}^\infty (1 - q^{n} )$.

Recall  that we have the following transformation formulas of theta-functions 
[{\bf Ch}]:
\begin{eqnarray}\begin{array}{l}
\theta (t+1, \tau) = -\theta (t,\tau), \qquad 
\theta ( t+ \tau, \tau)= - q^{-1/2} e^{- 2 \pi i t} \theta (t, \tau),\\
\theta_1 (t+1, \tau) = -\theta_1 (t,\tau), \qquad 
\theta_1 ( t+ \tau, \tau)=  q^{-1/2} e^{- 2 \pi i t} \theta_1 (t, \tau).
\end{array}\end{eqnarray}
and 
\begin{eqnarray}\begin{array}{l}
\theta ({t \over \tau}, - {1 \over \tau})= {1 \over i} \sqrt{\tau \over i}
 e^{\pi i t^2 \over \tau} \theta (t,\tau), \quad 
 \theta (t, \tau+1) = e^{ \pi i \over 4} \theta (t, \tau),\\
\theta_1 ({t \over \tau}, - {1 \over \tau})=  \sqrt{\tau \over i}
 e^{\pi i t^2 \over \tau} \theta_2 (t,\tau), \quad 
 \theta_1 (t, \tau+1) = e^{ \pi i \over 4} \theta_1 (t, \tau),\\
\theta_2 ({t \over \tau}, - {1 \over \tau})=  \sqrt{\tau \over i}
 e^{\pi i t^2 \over \tau} \theta_1 (t,\tau), \quad 
 \theta_2 (t, \tau+1) =  \theta_3 (t, \tau),\\
\theta_3 ({t \over \tau}, - {1 \over \tau})=  \sqrt{\tau \over i}
 e^{\pi i t^2 \over \tau} \theta_3 (t,\tau), \quad 
 \theta_3 (t, \tau+1) =  \theta_2 (t, \tau).
\end{array}\end{eqnarray}

Let $g= e^{2 \pi i t}\in S^1$ be a topological generator of $S^1$.
Let $\{M_{\alpha}\} $ be the fixed submanifolds of the circle action. 
Then $\pi: M_{\alpha}\to B$ be a submersion with fibre $X_{\alpha}$.
We have the following $S^1$-equivariant decomposition  of $TX$
\begin{eqnarray}
TX_{|M_{\alpha}} = N_1 \oplus \cdots \oplus N_h \oplus TX_{\alpha},
\end{eqnarray}
Here $N_{\gamma}$ is a complex vector bundle such that $g$ acts on it 
by $e^{2 \pi i m_{\gamma} t}$.
We denote the Chern roots of $N_{\gamma}$ by $2 \pi i x_{\gamma} ^j$, 
and the Chern roots of $TX_{\alpha} \otimes_{\sR } \bC$ by $\{ \pm 2 \pi i y'_j\}$. 
Let $\dim_{\mbox{\scriptsize \bf C}} N_\gamma = d(m_\gamma)$, and  
$\dim X_\alpha = 2 k_\alpha$.

Let 
\begin{eqnarray}
V_{|M_\alpha} = V_1 \oplus \cdots \oplus V_{l_0},
\end{eqnarray}
be the equivariant decomposition of $V$ restricted to $M_\alpha$. Assume that
$g$ acts on $V_v$ by $e^{2 \pi i n_v t}$, where some $n_v$ may be zero. 
We denote the Chern roots of $V_v$ by $2 \pi i u^j_v$. Let us write
$\dim_{\mbox{\scriptsize \bf R}} V_v = 2d(n_v)$.
By [{\bf Liu4}, \S 3.5],  the equivariant Chern character of $\psi (E,V)$ 
can be obtained as $q^{-m_\Lambda} C_{E, V} ( u+t, \tau)$, where 
\begin{eqnarray}
C_{E, V} ( u+t, \tau) = \chi_E (U+T,\tau),
 \end{eqnarray}
with  $ U+T= (u_1^j + n_1^j t, \cdots, 
u_{l_0}^j + n_{l_0}^j t)$.

For $g=e^{2 \pi i t}, t\in \bR$, and 
$\tau \in \bH$, $q= e^{ 2\pi i \tau}$, we let
\begin{eqnarray} 
F_{E,V}(t, \tau)=q^{m_\Lambda} \ch_g \Big ({\rm Ind} 
(D^X \bigotimes _{m=1}^\infty S_{q^m} (TX-\dim X) \otimes \psi (E, V)) \Big ). 
\end{eqnarray}

For $f(x)$ a holomorphic function, we denote by $f(y')(TX^g) = \Pi_j f(y'_j)$, 
the symmetric polynomial which gives characteristic class of $TX^g$,
 and similarly for $N_\gamma$. 
Using the family Atiyah-Bott-Segal-Singer Lefschetz fixed point formula
[{\bf LM}, Theorem 1.1], (1.14), we find for $t\in [0,1]\setminus \bQ$ 
\begin{eqnarray}  
\qquad F_{E,V}(t, \tau)= (2\pi i)^{-k} 
\sum_\alpha \pi_* \Big [ \theta' (0, \tau)^k
\Big ({2 \pi i y' \over \theta (y', \tau)}\Big )(TX^g) 
 {C_{E,V} ( u + t, \tau)
\over \Pi_{\gamma} \theta (x_{\gamma} + m_{\gamma} t, \tau)
(N_{\gamma})}\Big ].
\end{eqnarray}

Considered  as functions of $(t, \tau)$, we can obviously  extend 
$F_{E,V}( t, \tau)$ to meromorphic functions on $\bC \times \bH$ with values in
$H^*(B)$, 
and holomorphic in $\tau$. Theorem 1.2 is  equivalent to 
the statement that   $F_{E,V}(t, \tau)$ 
is independent of $t$. To prove it, we will proceed as in [{\bf Liu4}],
[{\bf LM}].

\begin{lemma} If $p_1(TX)_{S^1} = m p_1 (V)_{S^1}$, then for $a,b\in 2 \bZ$,
\begin{eqnarray}
F_{E,V}(t+ a\tau +b, \tau)= F_{E,V}(t, \tau).
\end{eqnarray}
\end{lemma}

$Proof$: By (1.15), for $a,b\in 2 \bZ$, $l\in \bZ$,  we have 
\begin{eqnarray}
\theta (x + l(t + a \tau + b), \tau) = e^{- \pi i (2l a x + 2 l^2 a t 
+ l^2 a^2 \tau)} \theta (x + lt , \tau).
\end{eqnarray}

Since $m p_1(V)_{S^1} = p_1(TX)_{S^1}$, we have
\begin{eqnarray}\begin{array}{l}
m  \Sigma_{v,j} (u_v^j + n_v t)^2 
= \Sigma_j (y'_j)^2 + \Sigma_{\gamma,j} (x_{\gamma}^j + m_{\gamma} t)^2.
\end{array}\end{eqnarray}
This implies the equalities:
\begin{eqnarray}\begin{array}{l}
m  \Sigma_{v,j} (u_v^j)^2=  \Sigma_j (y'_j)^2 
+ \Sigma_{\gamma,j} (x_{\gamma}^j)^2,\\
 m \Sigma_{v,j} n_v u_v^j = \Sigma_{\gamma,j} m_{\gamma} x_{\gamma}^j, \qquad 
\Sigma_{\gamma} m_{\gamma}^2 d(m_{\gamma}) = m \Sigma_v n_v^2 d(n_v).
\end{array}\end{eqnarray}
By using (1.11), (1.21),  (1.23), and (1.25),
 we get (1.22).\hfill $\blacksquare$\\

Now we will prove that $F_{E,V}(t, \tau)$ is holomorphic in $t$. 
Then, by Lemma 1.1, we get  the rigidity theorem.

To prove $F_{E,V}(t, \tau)$ is holomorphic in $t$, we will examine 
the modular transformation property of $F_{E,V}(t, \tau)$  
under the group $SL_2(\bZ)$.

Recall that for $g= \left ( \begin{array}{l} a \quad b\\
c \quad d
\end{array} \right ) \in SL_2 (\bZ)$,  we define its modular 
transformation on $\bC \times \bH$ by
\begin{eqnarray}\begin{array}{l}
\displaystyle{
g(t, \tau) = \left ( { t \over c \tau + d}, {a \tau + b \over c \tau + d}
\right ). }
\end{array}\end{eqnarray}

Obviously, the two generators of $SL_2(\bZ)$ are 
$ S = \left (  \begin{array}{l} 0 \quad -1 \\
1 \quad \quad 0
 \end{array} \right ), \
T = \left (  \begin{array}{l} 1 \quad 1 \\
0 \quad 1
 \end{array} \right ).$
They act on $\bC \times \bH$ in the following way:
\begin{eqnarray}
S(t,\tau) = \Big ({t \over \tau}, - {1 \over \tau}\Big ),\  
T(t, \tau) = (t, \tau +1).
\end{eqnarray}

Let $\Psi_{\tau}$ be the scaling homomorphism from $\Lambda (T^* B)$ into 
itself $: \beta \to \tau^{{1 \over 2} {\rm \scriptsize deg}\beta} \beta$.
If $ \alpha$ is a  differential form on $B$, we denote by $\{ \alpha\}^{(p)}$
the component of degree $p$ of $\alpha$.

\begin{lemma} For any $g= \left ( \begin{array}{l} a \quad b\\
c \quad d
\end{array} \right ) \in SL_2 (\bZ)$, we have 
\begin{eqnarray}
F_{E, V} (g(t, \tau)) = (c\tau +d)^k \Psi_{c\tau +d} F_{gE, V} (t, \tau),
\end{eqnarray}
where $gE= \Sigma_{\mu} a_{\mu} E_{\mu}$ is a finite complex linear 
combination of positive energy representations of $\widetilde{L} Spin(2l)$ 
of highest weight of level $m$, and we denote by 
\begin{eqnarray}\begin{array}{l}
\displaystyle{ F_{gE, V} (t, \tau) = (2 \pi i )^{-k} \theta' (0, \tau)^k
\sum_\mu \sum_{\alpha} a_{\mu} \pi_* \Big [ 
\Big ({2 \pi i y' \over \theta (y', \tau)}\Big )(TX^g) } \\
\displaystyle{\hspace*{40mm} 
{ C_{E_\mu,V} (u + t, \tau)\over \Pi_{\gamma} \theta (x_{\gamma} + m_{\gamma} t, \tau)(N_{\gamma})}\Big ].  }
\end{array}\end{eqnarray}
 the complex linear combination of the   equivariant 
Chern characters of the corresponding index bundles.
\end{lemma}

$Proof$: Set 
\begin{eqnarray}
F(t,\tau) = {\theta' (0, \tau)\over \theta (t,\tau)}.
\end{eqnarray}

By (1.16), we get
\begin{eqnarray}
F(g(t, \tau))= (c \tau + d) e ^{-c \pi i t^2/(c \tau + d)} 
F((c \tau + d) t ,\tau).
\end{eqnarray}
By Theorem 1.1, (1.19), it is easy to see that on $M_{\alpha}$,
\begin{eqnarray}
C_{E,V}(g(u+t, \tau))=e^{c m \pi i \Sigma_{v, j}
 (u_v^j + n_v t)^2/(c \tau + d)} C_{gE, V} (u+t, \tau).
\end{eqnarray}
with 
\begin{eqnarray}
C_{gE, V} (u+t, \tau)= \Sigma_{\mu} a_{\mu} C_{E_{\mu},V} (u + t, \tau).
\end{eqnarray}

By using (1.21), (1.31), (1.32), we get
\begin{eqnarray}\begin{array}{l}
\displaystyle{
F_{E, V} ({ t \over c \tau + d}, { a \tau + b \over c \tau + d}) 
= (2 \pi i)^{-k} \sum_\alpha \pi_*  \Big [\Big (2 \pi i y'
F(y', { a \tau + b \over c \tau + d} )\Big  )(TX^g)  }\\
\hspace*{25mm}\displaystyle{
 \Pi_{\gamma} \Big (F(x_{\gamma} + {m_{\gamma} t \over c \tau + d}, 
{ a \tau + b \over c \tau + d}) (N_{\gamma})\Big )
C_{E, V} (u + {t \over c \tau + d}, { a \tau + b \over c \tau + d}) \Big ] }\\
\displaystyle{
\hspace*{10mm} = (c \tau + d)^k (2 \pi i)^{-k} \sum_\alpha \pi_*  
\Big [\Big (2 \pi i y' 
F((c \tau + d)y', \tau) \Big  )(TX^g)  }\\
\hspace*{25mm} 
\displaystyle{\Pi_{\gamma} 
\Big (F((c \tau + d)x_{\gamma}+ m_{\gamma} t, \tau) (N_{\gamma})\Big )
C_{gE, V} ((c \tau + d)u + t, \tau)\Big ]  }
\end{array}\end{eqnarray}

By  (1.34), to prove (1.28),
  we only need prove the following equation for $p\in \bN$,
\begin{eqnarray}\begin{array}{l}
\displaystyle{
\Big \{\pi_*  \Big [\Big (2 \pi i y '
F((c \tau + d)y', \tau) \Big  )(TX^g) }\\
\hspace*{10mm}\displaystyle{
 \Pi_{\gamma} 
\Big (F((c \tau + d)x_{\gamma}+ m_{\gamma} t, \tau) (N_{\gamma})\Big )
C_{gE, V} ((c \tau + d)u + t, \tau)\Big ]  \Big \}^{(2p)}    }\\
\displaystyle{
= (c \tau + d) ^p \Big \{\pi_*  \Big [\Big (2 \pi i y' 
F(y', \tau) \Big  )(TX^g) }\\
\hspace*{25mm} \displaystyle{ \Pi_{\gamma} 
\Big (F(x_{\gamma}+ m_{\gamma} t, \tau) (N_{\gamma})\Big )
C_{gE, V} (u + t, \tau)\Big ]  \Big \}^{(2p)}    }.
\end{array}\end{eqnarray}
By looking at the degree $2(p+ k_{\alpha})$ part, that is the 
$(p+k_{\alpha})$-th homogeneous terms of the polynomials in $x$'s, $y'$'s and 
$u$'s, on both sides, we get (1.35). The proof of Lemma 1.2 is complete.\hfill $\blacksquare$\\

The following lemma is a generalization of [{\bf Liu4}, Lemma 2.3],
\begin{lemma} For any $g\in SL_2(\bZ)$, the function $F_{g E, V}(t, \tau)$ 
is holomorphic in $(t, \tau)$ for $(t, \tau) \in \bR \times \bH$.
\end{lemma}

$Proof$. Let $z= e^{2 \pi i t}$,
and let $N= {\rm max}_{\alpha, \gamma} |m_\gamma|$. 
Denote by $D_N\subset \bC ^2$ the domain 
\begin{eqnarray}\begin{array}{l}
|q|^{1/N} < |z| < |q|^{-1/N}, 0< |q| < 1.
\end{array}\end{eqnarray}
By (1.14), (1.19), (1.21)  and (1.29), we know that in $D_N$, 
$F_{g E, V}(t, \tau)$ has a convergent Laurent series expansion of the form 
\begin{eqnarray}
\sum_\mu a_\mu q^{m_{\Lambda_\mu}} \sum_{j=0}^\infty b^g_{j \mu} (z) q^j
\end{eqnarray}
Here $\{b^g_{j \mu} (z)\}$ are rational functions of $z$ with possible 
poles on the unit circle.

Now considered as a formal power series of $q$,
$$\bigotimes _{n=1}^{\infty} S_{q^n} (TX-\dim X) \otimes \Big ( \sum_\mu 
a_\mu q^{m_{\Lambda_\mu } } \psi (E_{\mu}, V) \Big ) = \sum_\mu a_{\mu} 
q^{m_{\Lambda_\mu } }\bigotimes _{j=0}^{\infty} V^g_{j, \mu} q^j$$
with  $V^g_{j, \mu}\in K_{S^1}(M)$. Note that the  terms in the above two 
sums correspond to each other. Now, we apply the family Lefschetz fixed 
point formula [{\bf LM}, Theorem 1.1] to each $V^g_{j, \mu}$, 
 for $t \in \bR \setminus \bQ$, we get
\begin{eqnarray}
b^g_{j \mu} (z) = \ch_z ({\rm Ind} (D \otimes V^g_{j, \mu})).
\end{eqnarray}
But by [{\bf S}, Proposition 2.2], we know that 
\begin{eqnarray}
K_{S^1}(B) \simeq K(B) \otimes R(S^1)
\end{eqnarray}
This implies that for  $t \in \bR \setminus \bQ, z= e^{ 2 \pi i t}$,
\begin{eqnarray}
b^g_{j \mu} (z) = \Sigma_{l= -N(j)}^{N(j)} a_{l,j}^{g, \mu} z^l.
\end{eqnarray}
for $N(j)$ some positive integer depending on $j$ 
and $a_{l,j}^{g, \mu}\in H^*(B)$.
Since both sides are analytic functions of $z$, this equality holds for 
any $z\in \bC$.

On the other hand, by multiplying $F_{g E, V}(t, \tau)$ by
$f(z) = \Pi_{\alpha, \gamma} (1- z^{m_\gamma})^{l' d(m_\gamma)}$ $(l'= \dim M)$, we get  holomorphic functions which have a convergent power series expansion of the form
$\Sigma_\mu a_\mu q^{m_{\Lambda_\mu}} \Sigma_{j=0}^\infty c_{j\mu}^g (z) q^j$,
 with $\{c^g_{j\mu}(z)\}$ polynomial functions in $D_N$. 
Comparing the above two expansions, one gets
\begin{eqnarray}
c^g_{j\mu}(z)= f(z) b^g_{j\mu}(z)
\end{eqnarray} 
for each $j$. So by the Weierstrass preparation theorem, we get $F_{gE,V}(t, \tau)$ is holomorphic in $D_N$.\hfill $\blacksquare$\\

{\em Proof of Theorem 1.1}: 
We will prove that  $F_{E,V}$ is holomorphic on 
$\bC \times \bH$, which implies the rigidity theorem we want to prove.

>From their expressions, we know the possible polar 
divisors of $F_{E,V}$ in $\bC \times \bH$ are of the form 
$t= {n \over l} (c\tau + d)$ 
with $n, c, d, l$ intergers and $(c,d)=1$ or $c=1$ and $d=0$. 

We can always find intergers $a, b$ such that $ad-bc = 1$, and consider 
the matrix 
$g=\left ( \begin{array}{l}d\quad -b\\
-c \quad a
\end{array} \right ) \in SL_2(\bZ)$. By (1.28),
\begin{eqnarray}\begin{array}{l}
\displaystyle{
\Psi_{(-c \tau + a)} F_{g E,V}(t,\tau) =  
(-c \tau + a)^{-k} F_{E,V} \Big ({t \over -c \tau + a},
 {d \tau -b \over -c \tau + a} \Big )}
\end{array}\end{eqnarray}

Now, if $t= {n \over l} (c\tau + d)$ is a polar divisor of $F_{E,V} (t,\tau)$, 
then one polar divisor of $F_{g E,V}(t,\tau)$ is given by
\begin{eqnarray}
{t \over -c \tau + a}= {n \over l}
 \Big ( c {d \tau -b \over -c \tau + a} + d \Big ),
\end{eqnarray}
which exactly gives $t = n/l$. This contradicts Lemma 1.3, and completes the proof of Theorem 1.1. \hfill $\blacksquare$

\subsection{ \normalsize Family vanishing theorems}

Recall that a (meromorphic) Jacobi form  of index $n$ and weight $l$ over 
$L\rtimes \Gamma$, where $L$ is an integral lattice  in the complex plane 
$\bC$ preserved by the modular subgroup $\Gamma \subset SL_2(\bZ)$,
 is a (meromorphic) function $F(t, \tau)$ on $\bC \times \bH$ such that
\begin{eqnarray}\begin{array}{l}
\displaystyle{F({t \over c \tau + d}, { a \tau + b \over c \tau +d}) 
= (c \tau + d)^l e^{2 \pi i n ( c t^2 / (c \tau + d))} F(t, \tau),}\\
\displaystyle{F(t + \lambda \tau + \mu, \tau) = 
e ^{- 2 \pi i n ( \lambda ^2 \tau + 2 \lambda t)} F(t, \tau),}
\end{array}\end{eqnarray}
where $(\lambda, \mu) \in L$, and $g = \left (\begin{array}{l} a \quad b\\
c \quad d
\end{array}  \right ) \in \Gamma$. If $F$ is holomorphic on $\bC \times \bH$, 
we say that $F$ is a holomorphic Jacobi form.

For $N\in \bN^*$, set
\begin{eqnarray}
\Gamma(N)= \left \{ g= \left ( \begin{array}{l} a \quad b\\
c \quad d  \end{array}\right ) \in SL_2 (\bZ)| 
g \equiv   \left ( \begin{array}{l} 1 \quad 0\\
0 \quad 1 \end{array} \right ) ({\rm mod} N) \right \}.
\end{eqnarray}

Recall that the equivariant cohomology group $H^*_{S^1} (M, \bZ)$ 
of $M$ is defined by
\begin{eqnarray}
H^*_{S^1} (M, \bZ)= H^*(M \times_{S^1} ES^1, \bZ).
\end{eqnarray}
where $ES^1$ is the usual  $S^1$-principal bundle over
the  classifying space $BS^1$ of $S^1$.
So $H^*_{S^1} (M, \bZ)$ is a module over $H^*(BS^1, \bZ)$ induced by the 
projection $\overline{\pi} : M\times _{S^1} ES^1\to BS^1$. 
Let $p_1(V)_{S^1}, p_1(TX)_{S^1} \in H^*_{S^1} (M, \bZ)$ be the equivariant
 first Pontrjagin classes of $V$ and $TX$ respectively.
Also recall that 
\begin{eqnarray}
H^*(BS^1, \bZ)= \bZ [[u]]
\end{eqnarray}
with $u$ a generator of degree $2$.

In this part, we suppose that there exists $n\in \bZ$ such that 
\begin{eqnarray}
m p_1(V)_{S^1}- p_1(TX)_{S^1} = n \cdot  \overline{\pi}^* u^2 
\quad {\rm  in} \quad H^*_{S^1} (M, \bZ)\otimes_{\sZ} \bQ .
\end{eqnarray}
As in [{\bf Liu4}], we call $n$ the anomaly to rigidity.

\begin{thm} Let $M,B, V$ and $ E$  be as in Theorem 1.2. Then for $p\in \bN$, 
  $\{F_{E, V}\}^{(2p)}$ is a holomorphic Jacobi form of index $n/2$ 
and weight $k+p$ over $(2 \bZ)^2 \rtimes  \Gamma(N(m))$. 
\end{thm}

Here $N(m)$ is an integer depending on the level $m$ and was given in [{\bf Kac}], and $m_{\Lambda}$ defined in (1.6).

$Proof$: Now, by (1.48), we get 
\begin{eqnarray}
m \Sigma_{v,j} (u_v^j + n_v t)^2 -\Big ( \Sigma_j (y'_j)^2 
+ \Sigma_{\gamma,j} (x_{\gamma}^j + m_{\gamma} t)^2\Big ) = n\cdot t^2.
\end{eqnarray}
This means
\begin{eqnarray}\begin{array}{l}
m \Sigma_v n_v^2 d(n_v)-\Sigma_{\gamma} m_{\gamma}^2 d(m_{\gamma}) =n, \quad 
m \Sigma_{v,j} n_v u_v^j = \Sigma_{\gamma,j} m_{\gamma} x_{\gamma}^j,\\
m \Sigma_{v,j} (u_v^j)^2 =
\Sigma_j (y'_j)^2 + \Sigma_{\gamma,j} (x_{\gamma}^j)^2.
\end{array}\end{eqnarray}

First by using (1.50), as in Lemma 1.1, for $(a,b)\in (2 \bZ)^2$, we get
\begin{eqnarray}
F_{E, V}(t + a \tau + b, \tau) = 
e^{- \pi i n (a^2 \tau + 2 a t) } F_{E,V}(t, \tau).
\end{eqnarray}

Second by a theorem of Kac, Peterson and Wakimoto [{\bf Kac}, Chapter 13],
there exists an integer $N(m)$ such that  for any 
$g=\left ( \begin{array}{l} a \quad b\\
c \quad d \end{array} \right )\in \Gamma(N(m))$, we have
\begin{eqnarray}
C_{E,V}(g(u+t, \tau)) = e^{m c \pi i \Sigma_{v, j} 
(u_v^j + n_v t)^2 /(c \tau + d)} C_{E, V} (u+ t, \tau).
\end{eqnarray}

Now, by using (1.21), (1.50) and  (1.52), as Lemma 1.2, we get
\begin{eqnarray}
F_{E,V}(g(t, \tau)) = (c \tau + d)^k e^{\pi i n c t^2/(c \tau + d)}
   \Psi_{c \tau + d} F_{E,V}(t, \tau).
\end{eqnarray}

As the same argument in the proof of Theorem 1.2 
(see also [{\bf Liu4}, Theorem 3.4]), we know $F_{E,V}$ is holomorphic on 
$\bC \times \bH$. By (1.51), (1.53), we get Theorem 1.3.\hfill $\blacksquare$\\

The following lemma was established in [{\bf EZ}, Theorem 1.2]:

\begin{lemma} Let $F$ be a holomorphic Jacobi form of index $m$ and weight $k$.
 Then for fixed $\tau$, $F(t,\tau)$, if not identically zero, has exactly 
 $2m$ zeros in any fundamental domain for the action of the lattice on $\bC$.
\end{lemma}
This tells us that there are no holomorphic Jacobi forms of negative index.
 Therefore, if $m<0$, $F$ must be identically zero. 
If $m=0$, it is easy to see  that $F$ must be independent of $t$. So we immediately get the following: 

\begin{cor} Let $M,B,V,E$ and $n$ be as in Theorem 1.3. If $n=0$, 
the equivariant Chern character of the index bundle  of
$$ {\rm Ind} \Big (D^X \bigotimes _{m=1}^\infty S_{q^m} (TX-\dim X) 
\otimes \psi (E, V)\Big ) $$
is independent of $g\in S^1$. If $n<0$, this equivariant Chern character 
is identically zero, in particular, the Chern character of this index bundle
is zero.
\end{cor}
 \hfill $\blacksquare$

\section{ \normalsize Family rigidity theorems for $spin^c$-manifolds}
\setcounter{equation}{0}

The purpose of this Section is to prove  a family version of  the rigidity 
theorem for $spin^c$-manifolds. 

This Section is organized as follows: 
In Section 2.1, we explain the equivariant family index theorem 
for $spin^c$-manifolds.
In Section 2.2, we state our main result, Theorem 2.2.
 In Section 2.3, we prove  Theorem 2.2.
In Section 2.4,  we prove a family version of  the rigidity and vanishing 
theorem for $spin^c$-manifolds of  [{\bf D1}]. 
A family vanishing theorem of Witten genus for $spin^c$-manifolds is also obtained.

\subsection{\normalsize Equivariant family index theorem 
for $spin^c$-manifolds}

By [{\bf LM}, Theorem 1.1], we have the equivariant family index theorem
for a family of  equivariant  elliptic operators. 
In fact, by the proof of [{\bf LM}], we know we have a local version of 
[{\bf LM}, Theorem 1.1] for the Dirac operator associated to 
 the Clifford module in the sense of  [{\bf BeGeV}, \S 3.3, \S 10.3].

Let $\pi: M\to B$ be a  fibration of compact manifolds with fiber $X$ with $\dim X= 2k$. We assume that the $S^1$ acts fiberwisely on $M$,
and $TX$ has an $S^1$-equivariant $spin^c$ structure.
Let $\Delta(TX)$ be the complex spinor bundle for $TX$ 
[{\bf LaM}, Definition D.9]. 
We denote $D^c$ the corresponding $spin^c$-Dirac operator on the fibre $X$
[{\bf LaM}, Appendix D].

Let $W$ be an $S^1$-equivariant complex vector bundle on $M$. 
Let $D^c \otimes W$ be the twisted
$spin^c$-Dirac operator on $\Delta(TX) \otimes W$.
Then ${\rm Ind} (D^c \otimes W)\in K_{S^1} (B)$.

Let $g= e^{2 \pi i t}\in S^1$ be a generator of the action group.
Let $\{M_{\alpha}\} $ be the fixed submanifolds of the circle action. 
Then $\pi: M_{\alpha}\to B$ be a submersion with fibre $X_{\alpha}$.
We have the following equivariant decomposition  of $TX$
\begin{eqnarray}
TX_{|M_{\alpha}} = N_1 \oplus \cdots \oplus N_h \oplus TX_{\alpha},
\end{eqnarray}
Here $N_{\gamma}$ is a complex vector bundle such that $g$ acts on it 
by $e^{2 \pi i m_{\gamma} t}$. So $TX_{\alpha}$ is naturally oriented.
We denote the Chern roots of $N_{\gamma}$ by $2 \pi i x_{\gamma} ^j$, 
and the Chern roots of $TX_{\alpha} \otimes_{\sR} \bC$ by $\{ \pm 2 \pi i y'_j\}$. 
Let $\dim_{\mbox{\scriptsize \bf C}} N_\gamma = d(m_\gamma)$, 
$\dim X_\alpha = 2 k_\alpha$.

We recall that the $spin^c$-structure on $TX$  induces an $S^1$-equivariant 
complex line bundle $L$ over $M$. 
Its equivariant Chern class $c_1(L)_{S^1}$
 will also be denoted by $c_1(TX)_{S^1}$. We denote the  Chern class $c_1(L)$ 
 of $L$ by $2 \pi i c_1$. Let $i_{\alpha}: M_{\alpha}\to M$ be the inclusion, 
and let $i^*_{\alpha}$ denote the induced homomorphism in equivariant
 cohomology.
If $g$ acts on $L$ on $M_\alpha$ by $e^{ 2 \pi i l_c t}$, we have
\begin{eqnarray}
i^*_{\alpha} c_1(TX)_{S^1} = 2 \pi i (c_1 + l_c t).
\end{eqnarray}

We denote $\pi_* : H^*(M^g) \to H^*(B)$ the intergration along the fibre $X^g$.
Now, we can reformulate the family Atiyah-Bott-Segal-Singer 
Lefschetz fixed point formula,  [{\bf LM}, Theorem 1.1] in this case,

\begin{thm} We have the following identity in $H^*(B)$
\begin{eqnarray}
\ch_g({\rm Ind} (D^c \otimes W))= \pi_* \left \{
{\widehat{A} (TX^g) \ch_g (W) 
e^{\pi i (c_1 + l_c t)}  \over  \Pi_{\gamma} \Big ( e^{ \pi i (x_{\gamma}
+m_{\gamma} t)}-
 e^{- \pi i (x_{\gamma}+m_{\gamma} t)} \Big )(N_{\gamma})} \right \}.
\end{eqnarray}
\end{thm}

\subsection{\normalsize Family rigidity for $spin^c$-manifolds}

In this part, we use the assumption  of  Section 2.1,
we also use the  notation of  Sections 1 and 2.1.

 For a vector bundle $E$ on $M$,
we denote by $\widetilde{E}$ the reduced vector bundle $E-\dim (E)$.

Let $W$ be an $S^1$-equivariant  complex vector bundle of rank $r$ over $M$.
 Let $L_W = \det (W)$ the determinant line bundle of $W$ on $M$.
Let $V$ be a dimension $2l$ real  vector bundle on $M$ with $S^1$-equivariant
$spin(2l)$ structure. Let $\Delta(V) = \Delta^+ (V) \oplus \Delta ^- (V) $
 be the spinor bundle of $V$.

Let $y= e^{ 2 \pi i \beta}$ be a complex number, and we define 
the following elements in $  K(M)[[q^{1/2}]]$:
\begin{eqnarray}\begin{array}{l}
\displaystyle{
\Theta_q (TX|W)_v=\bigotimes_{m=1}^{\infty}
S_{q^m} (\widetilde{TX } ) \otimes \Lambda_{{-1}} (W^*) 
\otimes \bigotimes_{n=1}^{\infty} 
\Lambda_{-q^n} (\widetilde{W \otimes_{\sR} \bC)}, }\\
\displaystyle{
\Theta_q^{\beta} (TX|W)_v = \bigotimes_{m=1}^{\infty}
S_{q^m} (\widetilde{TX } ) \otimes \Lambda_{-y^{-1}} (\widetilde{W}^*) 
\otimes \bigotimes_{n=1}^{\infty} \Lambda_{-yq^n} (\widetilde{W})
\otimes \Lambda_{-y^{-1}q^n} (\widetilde{W}^*).  }
\end{array}\end{eqnarray}
Let 
\begin{eqnarray}\begin{array}{l}
R_1(V)_v = \Delta (V)
\otimes \bigotimes_{n=1}^{\infty} \Lambda_{q^n} (\widetilde{V}),\\
R_2(V)_v =\bigotimes_{n=1}^{\infty} \Lambda_{-q^{n-{1 \over 2}} }
(\widetilde{V}),\\
R_3(V)_v =\bigotimes_{n=1}^{\infty} \Lambda_{q^{n-{1 \over 2}}}
(\widetilde{V}).
\end{array}\end{eqnarray}

For $g= e^{ 2 \pi i t}, t\in \bR$, $q=e^{2 \pi i \tau}, \tau \in \bH$, set 
\begin{eqnarray}\begin{array}{l}
F_1(t, \tau ) = 2^{-l} \ch_g  \Big ({\rm Ind} 
(D^c \otimes \Theta_q (TX|W)_v  \otimes R_1(V)_v) \Big ),\\
F^{\beta}_1(t, \tau ) = 2^{-l} \ch_g  \Big ({\rm Ind} 
(D^c \otimes \Theta_q^{\beta} (TX|W)_v  \otimes R_1(V)_v) \Big ).
\end{array}\end{eqnarray}
We consider $F_1(t, \tau ), F^{\beta}_1(t, \tau )$ as functions on 
$(t, \tau) \in \bC \times \bH$ with values in $H^*(B)$.

Recall that for $n \in \bN^*$, 
\begin{eqnarray}
\Gamma_1(n)= \left \{ g= \left ( \begin{array}{l} a \quad b\\
c \quad d  \end{array}\right ) \in SL_2 (\bZ)| 
g \equiv   \left ( \begin{array}{l} 1 \quad *\\
0 \quad 1 \end{array} \right ) ({\rm mod} n)
 \right \}.
\end{eqnarray}

\begin{thm} If  $p_1(V+W-TX)_{S^1} = n \cdot \overline{\pi}^* u^2$
 $(n\in \bZ )$, 
$c_1(W)_{S^1}= c_1(TX)_{S^1}$ in $H^*_{S^1} (M, \bZ) \otimes_{\sZ} \bQ$. Then

i).  If  $c_1(W)\equiv 0 \ {\rm mod} (N)$ $(N \in \bN, N> 1)$, 
then for $y = e^{2 \pi i \beta}$ an $N$th root of unity, and  $p\in \bN$,
$\{F^{\beta}_1(t, \tau )\}^{(2p)}$ is a holomorphic Jacobi form of index $n/2$
and weight $(k+p)$  over $(2N\bZ)^2 \rtimes  \Gamma_1(2N)$.

ii). For $p\in \bN$, $\{F_1(t, \tau )\}^{(2p)}$is a holomorphic Jacobi 
form of index $n/2$ and weight $k-r+p$ over $(2 \bZ)^2 \rtimes \Gamma_1(2)$. 
If $V=0$, the same holds for $\Gamma_1(2)$ replaced by $SL_2(\bZ)$.
\end{thm}

{\bf Remark:} If we replace the condition $c_1(W)_{S^1}= c_1(TX)_{S^1}$
by $\omega_2(TX)= \omega_2(W)$. Let $F_1(t, \tau )$, and let $F^{\beta}_1(t, \tau ) $
be the equivariant Chern character of index bundles of 
$D^c  \otimes (L_W \otimes L^{-1})^{1/2} \otimes  \Theta_q (TX|W)_v 
 \otimes R_1(V)_v$, $D^c \otimes (L_W \otimes L^{-1})^{1/2} 
\otimes \Theta_q^{\beta} (TX|W)_v  \otimes R_1(V)_v$,
 we still have Theorem 2.2.  In fact, we only need to take tensor product of 
$(L_W \otimes L^{-1})^{1/2}$ with the  corresponding operators 
in the proof of Theorem 2.2. 
If $V=0$, this result also generalizes [{\bf Liu2}, Theorem B] to the  family case.

{\bf Remark:} If we replace the condition $c_1(W)_{S^1}= c_1(TX)_{S^1}$
by $c_1(W)= c_1(TX)$ in $H^*(M, \bQ)$, as Dessai remarked in [{\bf D}, Lemma 3.4], then there exists
$m\in \bZ$ such that $c_1(W)_{S^1}- c_1(TX)_{S^1}= m \overline{\pi}^* u$. 
Now, for the functions $e^{\pi i m t} F_1(t, \tau )$, 
$e^{\pi i m t} F^{\beta}_1(t, \tau )$, we still have the result 
of Theorem 2.2. In fact this only multiplies all the functions in the proof of 
Theorem 2.2 by $e^{\pi i m t}$.

\subsection{\normalsize Proof of Theorem 2.2}

Let
\begin{eqnarray}
V_{|M_\alpha} = V_1 \oplus \cdots \oplus V_{l_0},
\end{eqnarray}
be the equivariant decomposition of $V$ restricted to $M_\alpha$. Assume that
$g$ acts on $V_v$ by $e^{2 \pi i n_v t}$, where some $n_v$ may be zero. 
We denote the Chern roots of $V_v$ by $2 \pi i u^j_v$. Write
$\dim_{\mbox{\scriptsize \bf R}} V_v = 2 d(n_v)$.
Similarly, let
\begin{eqnarray}
W_{|M_\alpha} = W_1 \oplus \cdots \oplus W_{r_0},
\end{eqnarray}
be the equivariant decomposition of $W$ restricted to $M_\alpha$. Assume that
$g$ acts on $W_\mu$ by $e^{2 \pi i r_\mu t}$, where some $r_\mu$ may be zero. 
We denote the Chern roots of $W_\mu$ by $2 \pi i \omega^j_\mu$. Write
$\dim_{\mbox{\scriptsize \bf C}} W_\mu = d(r_ \mu)$.

First note that the condition $c_1(W)_{S^1}= c_1(TX)_{S^1}$ means that
\begin{eqnarray}
\Sigma_{\mu, j}(\omega^j_\mu+ r_\mu t) = c_1 + l_ct.
\end{eqnarray}

We take $ \beta = 1/N$. By 
applying the family Atiyah-Bott-Segal-Singer Lefschetz fixed point formula, 
Theorem 2.1, and using (2.10), 
for $g= e^{ 2 \pi i t}, t\in \bR \setminus \bQ$, we get
\begin{eqnarray}\begin{array}{l}
\displaystyle{F_1(t, \tau ) =  (2 \pi i)^{r-k} 
{ \theta'(0, \tau)^{k-r} \over \theta_1(0, \tau)^l}
\sum_{\alpha} \pi_* \Big [ 
\Big ({2 \pi i y' \over \theta (y', \tau)}\Big )(TX^g) }\\
\hspace*{20mm} \displaystyle{ {\Pi_v  \theta_1(u_v + n_v t, \tau) (V_v) 
\over \Pi_{\gamma}\theta (x_{\gamma} + m_{\gamma} t, \tau)
(N_{\gamma})} \Pi_\mu  \theta(\omega _\mu + r_\mu t, \tau) 
(W_\mu) \Big ],}\\
\displaystyle{F^{\beta}_1 (t, \tau ) =  (2 \pi i)^{-k} 
{ \theta'(0, \tau)^k \over \theta_1(0, \tau)^l \theta(\beta, \tau)^r}
\sum_{\alpha} \pi_* \Big [ 
\Big ({2 \pi i y' \over \theta (y', \tau)}\Big )(TX^g) }\\
\hspace*{20mm} \displaystyle{ {\Pi_v  \theta_1(u_v + n_v t, \tau) (V_v) 
\over \Pi_{\gamma}\theta (x_{\gamma} + m_{\gamma} t, \tau)
(N_{\gamma})} \Pi_\mu  \theta(\omega _\mu + r_\mu t+\beta, \tau) 
(W_\mu) \Big ].}
\end{array}\end{eqnarray}

In the following, we  will consider $F_1(t, \tau )$, $ F^{\beta}_1 (t, \tau )$
as meromorphic functions on  $(t, \tau) \in  \bC \times \bH$
 with values in $H^*(B)$.

\begin{lemma} If  $p_1(V+W-TX)_{S^1} = n \cdot \overline{\pi}^*( u^2)$ 
in $H^*_{S^1} (M, \bZ) \otimes_{\sZ} \bQ$ for some integer $n$,

i) For $a, b \in 2  \bZ$, 
 \begin{eqnarray}
F_1(t+ a \tau + b, \tau ) = 
e ^{-  \pi i n ( a^2 \tau + 2 a t)} F_1(t, \tau ).
\end{eqnarray}

ii) For $a, b \in 2 N \bZ$, 
\begin{eqnarray}
F^{\beta}_1(t+ a \tau + b, \tau ) = 
e ^{-  \pi i n ( a^2 \tau + 2 a t)} F^{\beta}_1(t, \tau ).
\end{eqnarray}
\end{lemma}

$Proof$:  Since  $p_1(V+W-TX)_{S^1} = n \cdot \overline{\pi}^* u^2$, we have
\begin{eqnarray}\begin{array}{l}
\Sigma_{v,j} (u_v^j + n_v t)^2+ \Sigma_{\mu, j}(\omega_\mu^j + r_\mu t)^2 
-( \Sigma_j (y'_j)^2 + \Sigma_{\gamma,j} (x_{\gamma}^j + m_{\gamma} t)^2)
= n t^2. 
\end{array}\end{eqnarray}

This implies the equalities:
\begin{eqnarray}\begin{array}{l}
\Sigma_{v,j} (u_v^j)^2 +  \Sigma_{\mu, j} (\omega_\mu^j)^2 =  \Sigma_j (y'_j)^2
+ \Sigma_{\gamma,j} (x_{\gamma}^j)^2,\\
 \Sigma_{v,j} n_v u_v^j+ \Sigma_{\mu, j} \omega_\mu^j r_{\mu}
 = \Sigma_{\gamma,j} m_{\gamma} x_{\gamma}^j,\\
\Sigma_v n_v^2 d(n_v) + \Sigma_{\mu}r_{\mu}^2 d(r_{\mu})
- \Sigma_{\gamma} m_{\gamma}^2 d(m_{\gamma}) =n.
\end{array}\end{eqnarray}

By (1.15), for $\theta_v=\theta, \theta_1$; $a,b\in 2 \bZ$, $l\in \bZ$, we have
\begin{eqnarray}
\theta_v (x + l(t + a \tau + b), \tau) = e^{- \pi i (2l a x + 2 l^2 a t 
+ l^2 a^2 \tau)} \theta_v (x + lt , \tau).
\end{eqnarray}

Let $F_{1,\alpha}$, $F^{\beta}_{1,\alpha}$ be the contribution 
of $M_{\alpha}$ to $F_1(t, \tau )$,  $F^{\beta}_1(t, \tau )$.
By using (2.11), (2.15), (2.16), we get for $a,b\in 2 \bZ$, 
\begin{eqnarray}\begin{array}{l}
F_{1,\alpha}(t + a \tau + b, \tau) = 
e^{ - \pi i n (a^2 \tau + 2 a t)} F_{1,\alpha} (t, \tau),\\
F^{\beta}_{1,\alpha}(t + a \tau + b, \tau) = y ^{- \Sigma_{\mu} r_{\mu} a}
e^{ - \pi i n (a^2 \tau + 2 a t)} F^{\beta}_{1,\alpha} (t, \tau).
\end{array}\end{eqnarray}
Since by the assumption, $y^N=1$, we get Lemma 2.1.
\hfill $\blacksquare$\\

For $A= \left ( \begin{array}{l} a \quad b\\
c \quad d
\end{array} \right ) \in SL_2 (\bZ)$, we define 
\begin{eqnarray}
\varepsilon_A =\left \{ \begin{array}{l} 1, \quad { \rm if} \quad (c,d)\equiv (0,1) \quad  ({\rm mod} 2) ,\\
2, \quad { \rm if} \quad (c,d)\equiv (1,0) \quad  ({\rm mod} 2) ,\\
3, \quad { \rm if} \quad (c,d)\equiv (1,1)\quad   ({\rm mod} 2) 
\end{array} \right.
\end{eqnarray}

For $g=e^{2 \pi i t}, t\in \bR$,
$A= \left ( \begin{array}{l} a \quad b\\
c \quad d
\end{array} \right ) \in SL_2 (\bZ)$, $j= 1,2,3$, we let
\begin{eqnarray}\begin{array}{l}
F_j (t, \tau )=\varepsilon (j)  \ch_g  \Big ({\rm Ind} 
(D^c \otimes \Theta_q (TX|W)_v  \otimes R_j(V)_v) \Big ),\\
F^{\beta}_j (t, \tau )^A =\varepsilon (j)  \ch_g  \Big ({\rm Ind} 
(D^c \otimes L_W ^{c \beta}\otimes \Theta_q^{(c \tau + d)\beta} (TX|W)_v  
\otimes R_j(V)_v) \Big ).
\end{array}\end{eqnarray}
with $\varepsilon (j)= 2^{-l}$ for $j=1$; $1$ for $j=2,3$.

By applying the family Atiyah-Bott-Segal-Singer Lefschetz fixed point formula, 
Theorem 2.1, and using (2.10), 
for $g= e^{ 2 \pi i t}, t\in \bR \setminus \bQ$, $j = 1,2,3$,  we get
\begin{eqnarray}  \qquad \begin{array}{l}
\displaystyle{F_j(t, \tau ) =  (2 \pi i)^{r-k} 
{ \theta'(0, \tau)^{k-r} \over \theta_j(0, \tau)^l}
\sum_{\alpha} \pi_* \Big [ 
\Big ({2 \pi i y' \over \theta (y', \tau)}\Big )(TX^g) }\\
\hspace*{20mm} \displaystyle{ {\Pi_v  \theta_j(u_v + n_v t, \tau) (V_v) 
\over \Pi_{\gamma}\theta (x_{\gamma} + m_{\gamma} t, \tau)
(N_{\gamma})} \Pi_\mu  \theta(\omega _\mu + r_\mu t, \tau) 
(W_\mu) \Big ],}\\
\displaystyle{F^{\beta}_j (t, \tau )^A =  (2 \pi i)^{-k} 
{\theta'(0, \tau)^k 
\over \theta_j(0, \tau)^l \theta((c \tau + d)\beta, \tau)^r}
\sum_{\alpha} \pi_* \Big [ 
\Big ({2 \pi i y' \over \theta (y', \tau)}\Big )(TX^g) }\\
\hspace*{5mm} \displaystyle{ {\Pi_v  \theta_j (u_v + n_v t, \tau) (V_v) 
\over \Pi_{\gamma}\theta (x_{\gamma} + m_{\gamma} t, \tau)
(N_{\gamma})} \Pi_\mu  \Big ( e^{2 \pi i  c\beta( \omega _\mu + r_\mu t)}   
\theta(\omega _\mu + r_\mu t+(c \tau + d)\beta, \tau) \Big )
(W_\mu) \Big ].}
\end{array}\end{eqnarray}

{\bf Remark:} In fact, to define a $S^1$-action on $L_W^{c \beta}$, 
we must replace the $S^1$- action by its $N$-fold action. Here by abusing notation, we still say an $S^1$-action without causing any confusion.

\begin{lemma} If $p_1(V+W-TX)_{S^1} = n \cdot \overline{\pi}^*( u^2)$,
under the action $A= \left ( \begin{array}{l} a \quad b\\
c \quad d \end{array} \right ) \in SL_2 (\bZ)$, we have 
\begin{eqnarray}\begin{array}{l}
F_1(A(t, \tau)) = (c\tau +d)^{k- r} 
e^{\pi i n c t^2/(c \tau + d)} \Psi_{c \tau +d} F_{\varepsilon_A} (t, \tau),\\
F^{\beta}_1(A(t, \tau)) = (c\tau +d)^k e^{\pi i n c t^2/(c \tau + d)}
 \Psi_{c \tau +d} F^{\beta}_{\varepsilon_A}(t, \tau)^A,\\
\end{array}\end{eqnarray}
\end{lemma}

$Proof$: As the equation for $F_1$ in (2.21) is very easy, 
we leave it  to the interested reader. 
Here, we only prove (2.21) for $F^\beta_1$.

 By (1.15), (1.16), ${\theta_1(t, \tau) \over \theta_1(0, \tau)}$
is a Jacobi form of index $1/2$ and weight $0$ over $(2 \bZ)^2 \rtimes \Gamma_1(2)$. This explains the index $\varepsilon_A$ in the following equation. 
 By (1.16),  we get
\begin{eqnarray}\begin{array}{l}
\displaystyle{
{\theta'(0, {a \tau + b \over c \tau + d}) \over 
\theta ({t \over c \tau + d}, {a \tau + b \over c \tau + d})} 
= (c \tau + d) e^{- \pi i { c t^2 \over c \tau + d}} 
{\theta'(0, \tau) \over \theta (t, \tau)},   }\\
\displaystyle{
{\theta({t \over c \tau + d}, {a \tau + b \over c \tau + d})
\over \theta( \beta , {a \tau + b \over c \tau + d})} =
e^{\pi i c ( {t^2 \over c \tau + d} - \beta^2 (c \tau + d))} 
{\theta(t, \tau) \over \theta(\beta (c \tau+ d), \tau)},   }\\
\displaystyle{
{\theta_1({t \over c \tau + d}, {a \tau + b \over c \tau + d})
\over \theta_1( 0, {a \tau + b \over c \tau + d})} =
e^{\pi i c {t^2 \over c \tau + d}} 
{\theta_{\varepsilon_A}(t, \tau) \over \theta_{\varepsilon_A} (0, \tau)}.}
\end{array}\end{eqnarray}
By (2.11), (2.22), we get 
\begin{eqnarray}\qquad \begin{array}{l}
\displaystyle{
F^{\beta}_1({ t \over c \tau + d}, { a \tau + b \over c \tau + d}) 
= (2 \pi i)^{-k} \sum_\alpha \pi_*  \Big [\Big (2 \pi i y'
{ \theta'(0,{ a \tau + b \over c \tau + d} ) 
\over \theta (y',{ a \tau + b \over c \tau + d} )}\Big )(TX^g)  }\\
 \displaystyle{\hspace*{20mm}
 \Pi_{\gamma}\Big ( {\theta'(0, { a \tau + b \over c \tau + d}) \over 
\theta (x_{\gamma} + m_{\gamma}  {t \over c \tau + d},
 { a \tau + b \over c \tau + d})}\Big ) (N_{\gamma}) 
\Pi_v \Big ( { \theta_1(u_v + n_v {t \over c \tau + d}, 
{ a \tau + b \over c \tau + d}) 
\over  \theta_1(0, { a \tau + b \over c \tau + d})}\Big ) (V_v) }\\
  \displaystyle{\hspace*{20mm}
\Pi_\mu { \theta(\omega _\mu + r_\mu  {t \over c \tau + d} +\beta, 
{ a \tau + b \over c \tau + d}) 
\over \theta(\beta, { a \tau + b \over c \tau + d})}(W_\mu) \Big ] }\\
\displaystyle{
= (2 \pi i)^{-k}  (c\tau +d)^k e^{\pi i n c t^2/(c \tau + d)}
\sum_\alpha \pi_*  \Big [\Big (2 \pi i y' { \theta'(0,\tau) 
\over \theta ( (c \tau + d) y', \tau)}\Big )(TX^g)  }\\
 \displaystyle{\hspace*{20mm}
 \Pi_{\gamma}\Big (  {\theta'(0,  \tau ) \over 
\theta ((c \tau + d) x_{\gamma} + m_{\gamma} t,\tau )}\Big )(N_{\gamma}) 
\Pi_v\Big ( { \theta_{\varepsilon_A}((c \tau + d) u_v + n_v t, \tau ) 
\over  \theta_{\varepsilon_A}(0, \tau )}\Big ) (V_v) }\\
  \displaystyle{\hspace*{20mm}
\Pi_\mu \Big [e^{2 \pi i c ((c \tau +d )\omega _\mu + r_\mu t) \beta}
  {\theta((c \tau + d) \omega _\mu + r_\mu t +(c \tau + d) \beta, \tau ) 
\over \theta((c \tau + d) \beta,\tau )} \Big ](W_\mu) \Big ] }.
\end{array}\end{eqnarray}
To prove (2.23), we will  prove 
\begin{eqnarray} \qquad \begin{array}{l}
 \displaystyle{ 
\Big \{\pi_*  \Big [\Big ( {2 \pi i y'
\over \theta ( (c \tau + d) y', \tau)}\Big )(TX^g)  
 {\Pi_v \theta_{\varepsilon_A}((c \tau + d) u_v + n_v t, \tau ) (V_v) \over 
\Pi_{\gamma}\theta ((c \tau + d) x_{\gamma} + m_{\gamma} t,\tau )(N_{\gamma})} 
}\\
  \displaystyle{\hspace*{10mm}
\Pi_\mu \Big [e^{2 \pi i c \beta ((c \tau +d )\omega _\mu+ r_\mu t)} 
\theta \Big ((c \tau + d)\omega_\mu + r_\mu t +(c \tau + d)\beta, \tau \Big )
 \Big ]  (W_\mu) \Big ] \Big \}^{(2p)}    }\\
\displaystyle{ 
= (c \tau + d)^{p}  \Big \{\pi_* \Big [ 
\Big ({2 \pi i y' \over \theta (y', \tau)}\Big )(TX^g) 
{\Pi_v  \theta_{\varepsilon_A}(u_v + n_v t, \tau) (V_v) 
\over \Pi_{\gamma}\theta (x_{\gamma} + m_{\gamma} t, \tau)(N_{\gamma})}  }\\
\hspace*{20mm} \displaystyle{ 
 \Pi_\mu  \Big [e^{2 \pi i c \beta (\omega _\mu+ r_\mu t)} 
\theta \Big (\omega _\mu + r_\mu t+(c \tau + d)\beta, \tau \Big ) \Big ] 
(W_\mu) \Big ]\Big \}^{(2p)}  .}
\end{array}\end{eqnarray}
  By looking at the degree $2(p+k_\alpha)$ part, that is the $(p+k_\alpha)$-th 
homogeneous terms of the polynomials in $x$'s,  $y'$'s, $u$'s and $\omega$'s   
 on both sides, we immediately get (2.24).

The proof of Lemma 2.2 is complete.\hfill $\blacksquare$\\

Since $F_j(t, \tau )$,  $F^{\beta}_j(t, \tau )^A$ $(j=1,2,3)$ are 
the equivariant Chern characters of the index bundles of some elliptic operators, the same proof as that of Lemma 1.3 gives the following 

\begin{lemma} 
i) $F_j(t, \tau )$ $(j=1,2,3)$ is holomorphic in $(t, \tau)\in \bR \times \bH$.

ii) If $c_1(W)\equiv 0 ({\rm mod} N)$, then for $A\in SL_2(\bZ)$, $j=1,2,3$, 
$F^{\beta}_j(t, \tau )^A$ is holomorphic in $(t, \tau)\in \bR \times \bH$.
\end{lemma}

This is the only essential place where we need the topological condition 
$c_1(W)\equiv 0 ({\rm mod} N)$ which insures the existence of $L^{c \beta}_W$,
therefore the holomorphicity of $F^{\beta}_1(t, \tau )^A$ for $t \in \bR$. 

{\em Proof of Theorem 2.2}: Now, if $A= \left (\begin{array}{l} a \quad b\\
c \quad d \end{array}  \right ) \in \Gamma_1(2N)$, by (1.15), (2.20),  we get 
\begin{eqnarray}
F^{\beta}_{\varepsilon_A}(t,\tau)^A = F^{\beta}_1(t, \tau).
\end{eqnarray}

By using the above three Lemmas, and proceeding as in 
the proof of Theorem 1.1, we know that $F^{\beta}_1(t, \tau )$
is holomorphic in $(t, \tau) \in \bC \times \bH$. 

>From Lemmas 2.1, 2.2, (2.25), we get Theorem 2.2.\hfill $\blacksquare$\\

\subsection{\normalsize  Family rigidity and vanishing 
theorems for $spin^c$-manifolds}

From Lemma 1.4 and Theorem 2.2, we get the following family rigidity and 
vanishing theorems for $spin^c$-manifolds.

\begin{thm} Let $M, B, W, V$ as in Theorem 2.2. 
If  $p_1(V+W-TX)_{S^1} = n \cdot \overline{\pi}^* u^2$ $(n\in \bZ)$
and $c_1(W)_{S^1}= c_1(TX)_{S^1}$ in  $H^*_{S^1} (M, \bZ) \otimes_{\sZ} \bQ$.

i). If $n=0$, then $D^c \otimes \Theta_q (TX|W)_v  \otimes R_1(V)_v$ 
is rigid. If, in addition, $c_1(W)$ is divisible by an integer $N\geq 2$, 
then $D^c \otimes \Theta_q^{\beta} (TX|W)_v  \otimes R_1(V)_v$ is rigid for 
$y= e^{ 2 \pi i \beta}$ an $N$th root of unity.

ii). If $n<0$, then  the equivariant Chern character of the index bundle 
$D^c \otimes \Theta_q (TX|W)_v$ $\otimes R_1(V)_v $ vanishes identically, 
in particular, the Chern character of this index bundle is zero.  
If, in addition, $c_1(W)$ is divisible by an integer $N\geq 2$,
then the equivariant Chern character of the index bundle
 $D^c \otimes \Theta_q^{\beta}(TX|W)_v  \otimes R_1(V)_v$ vanishes identically
for $y= e^{ 2 \pi i \beta}$ an $N$th root of unity,
 in particular, the Chern character of this index bundle is zero.
\end{thm}

The following family vanishing Theorem generalizes
[{\bf LM}, Theorem 3.2] to family $spin^c$-manifolds.

\begin{thm} Let $\pi: M \to B$ be a fibration of compact connected  manifolds 
with compact fibre $X$, and $S^1$ acts fiberwisely and non-trivially on $M$. 
We suppose $TX$ has a $S^1$-equivariant 
$spin^c$ structure. If $c_1(TX)=0$ in $H^*(M, \bQ)$, 
and if $p_1(TX)_{S^1} =- n \cdot \overline{\pi}^* u^2$ in $H^*_{S^1}(M, \bZ)\otimes_{\sZ} \bQ$ for some integer $n$, 
then the equivariant Chern character of the index bundle,
 especially the Chern character of the index bundle of 
$D^c \otimes \otimes _{m=1}^\infty S_{q^m} (\widetilde{TX})$ is zero. 
\end{thm}

{\bf Remark}: Note that the condition $c_1(TX)=0$ in $H^*(M, \bQ)$ does not mean the $Spin^c$ structure is spin. This only insures that there exists $m \in \bZ$, 
such that $c_1(TX)_{S^1} = m \overline{\pi}^* u$. So in fact, 
the difference between [{\bf LM}, Theorem 3.2] and Theorem 2.4 are quite subtle.

As pointed out by Dessai [{\bf D2}, \S 3], when the $S^1$-action is induced 
from an $S^3$ or nice Pin(2) action on $M$ (In fact the $S^3$ and $Pin(2)$ action need not act fiberwisely on $M$), the condition
 $p_1(TX)_{S^1} =- n \cdot \overline{\pi}^* u^2$ in $H_{S^1}(M, \bZ)\otimes_{\sZ} \bQ$   is also equivalent to $p_1(TX)=0$ in $H^*(M, \bQ)$.

In [HL], some related result was proved for foliations. 
\\

{\em Proof of Theorem 2.4}: We only need to put  $W=V=0$ in Theorem 2.2. In fact, by (2.15), we know 
\begin{eqnarray}
\Sigma_j m^2_j d(m_j)=-n.
 \end{eqnarray}
So  the case $n>0$ can never happen. If $n=0$, then all the exponents 
$\{m_j\}$ are zero, so the $S^1$-action can not  have a fixed point. 
By (2.11), we get the result.
For $n<0$, by Remark in Section 2.2 and Theorem 2.3, we get  the result.
\hfill $\blacksquare$\\

\centerline{------------------------}
\vskip 6mm

Kefeng LIU,
Department of Mathematics, Stanford University, Stanford, CA 94305, USA.

{\em E-mail address}: kefeng@math.stanford.edu

\vskip 6mm

Xiaonan MA,
Humboldt-Universitat zu Berlin,
Institut f\"ur Mathematik,unter den Linden 6,
D-10099 Berlin, Germany.

{\em E-mail address}: xiaonan@mathematik.hu-berlin.de

\end{document}